\documentclass[12pt,oneside,english]{amsart}
\usepackage[latin1]{inputenc}
\usepackage{amssymb}

\makeatletter
\newtheorem{theorem}{Theorem}
\newtheorem{lemma}{Lemma}
\newtheorem{example}{Example}
\newtheorem{definition}{Definition}

\newtheorem{remark}{Remark}
\newtheorem{conjecture}{Conjecture}

\numberwithin{equation}{section}
\usepackage{babel}

\makeatother
\begin{document}
\baselineskip=17pt

\title[Spectrum of permanent's values in $\Lambda_n^3$ and $\Lambda_n(\alpha,\beta,\gamma)$]{Spectrum of permanent's values and its extremal magnitudes in $\Lambda_n^3$ and $\Lambda_n(\alpha,\beta,\gamma)$}

\author{Vladimir Shevelev}
\address{Department of Mathematics \\Ben-Gurion University of the
 Negev\\Beer-Sheva 84105, Israel. e-mail:shevelev@bgu.ac.il}


\begin{abstract}
Let $\Lambda_n^k$ denote the class of $(0,1)$ square matrices containing in each row and in each column exactly $k$ 1's. The minimal value of $k,$ for which the behavior of the permanent in $\Lambda_n^k$ is not quite studied, is $k=3.$ We give a simple algorithm for calculation upper magnitudes of permanent in $\Lambda_n^3$ and consider some extremal problems in a generalized class $\Lambda_n(\alpha,\beta,\gamma),$ the matrices of which contain in each row and in each column nonzero elements $\alpha,\beta,\gamma$ and $n-3$ zeros.
\end{abstract}

\maketitle

\section{Introduction }
Let $\Lambda_n^k$ denote the class of $(0,1)$ square matrices containing in each row and in each column exactly $k$ 1's. If $A\in \Lambda_n^3,$ then matrix $ k^{-1}A$ is doubly stochastic. Therefore, $\Lambda_n^k$-matrices are also called \slshape doubly stochastic (0,1)-matrices \upshape \enskip(cf. \cite{11}). Furthermore, for a given real or complex numbers $\alpha,\beta,\ldots\gamma,$ denote $\Lambda_n(\alpha,\beta,\ldots\gamma)$ the class of square matrices containing every number from $\{\alpha,\beta,\ldots\gamma\}$ exactly one time in each row and in each column, such that the other elements are 0's.
\begin{definition}
We call \upshape $p$-spectrum \slshape in $\Lambda_n^k$ (denoting it $ps[\Lambda_n^k])$\enskip the set of all the values which are taken by the permanent in $\Lambda_n^k.$
\end{definition}
Note that $p$-spectrum in $\Lambda_n^1$ trivially is $\{1\}.$ It is known (cf. Tarakanov \cite{25}) that
 $$ps[\Lambda_n^2]=\{2,2^2,2^3,...,2^{\lfloor\frac{n} {2}\rfloor}\}.$$
  But, for $k\geq3,$  $p$-spectrum of $\Lambda_n^k,$ generally speaking, is unknown.
Greenstein (cf. \cite{11}, point 8.4, Problem 3) put the problem of describing the p-spectrum in $\Lambda_n^3.$
In this paper we find $p$-spectrum on symmetric matrices in $\Lambda_n^3$ with ones on the main diagonal and give an algorithm for calculation upper values of $p$-spectrum in $\Lambda_n^3.$ We also obtain several results for a generalized class $\Lambda_n(\alpha,\beta,\gamma)$ with real nonzero numbers $\alpha,\beta,\gamma.$ Some results of the present paper were announced by the author in \cite{21}.
\section{What is known about $\Lambda_n^3?$}
1) Explicit formula for $|\Lambda_n^3|$ (cf. Stanley \cite{24}, Ch.1)
 \begin{equation}\label{2.1}
 |\Lambda_n^3|=6^{-n}\sum_{k_1+k_2+k_3=n,\enskip k_i\geq0}\frac{(-1)^{k_2}n!^2(k_2+3k_3)!2^{k_1}3^{k_2}}{k_1!k_2!k_3!^26^{k_3}}.
 \end{equation}
\indent 2) Asymptotic formula for $|\Lambda_n^3|$ (cf. O'Neil \cite{12})
\begin{equation}\label{2.2}
|\Lambda_n^3|=\frac{(3n)!}{(36)^n}e^{-2}(1+O(n^{-1+\varepsilon})),
\end{equation}
where $\varepsilon>0$ is arbitrary small for sufficiently large $n.$ \newline
\indent In addition, note that $|\Lambda_n(\alpha,\beta,\gamma)|$ with different $\alpha,\beta,\gamma$ is, evidently, the number of $3$-rowed Latin rectangles of length $n$ such that
\begin{equation}\label{2.3}
|\Lambda_n(\alpha,\beta,\gamma)|=n!K_n,
\end{equation}
where $K_n$ is the number of reduced $3$-rowed Latin rectangles with the first row $\{1,2,...,n\}.$ It is known (Riordan \cite{13}, pp. 204-210) that
\begin{equation}\label{2.4}
K_n=\sum_{k=0}^{\lfloor\frac{n}{2}\rfloor}\binom{n} {k}D_{n-k}D_kU_{n-2k},
\end{equation}
where $D_n$ is subfactorial.
\begin{equation}\label{2.5}
D_0=1,\enskip D_n=nD_{n-1}+(-1)^n, \enskip n\geq1,
\end{equation}
$\{U_n\}$ is sequence of Lucas numbers of the M\'{e}nage problem which is defined by Cayley recursion (cf. \cite{13}, p. 201)
$$U_0=1,\enskip U_1=-1,\enskip U_2=0,$$
\begin{equation}\label{2.6}
\enskip U_n=nU_{n-1}+\frac {n}{n-2}U_{n-2}+4\frac{(-1)^{n}} {n-2},\enskip n\geq3
\end{equation}
(see \cite{23}, sequences A102761, A000186).\newline
\indent Denote, furthermore, $\overline{\Lambda}_n^3$ the set of matrices in  $\Lambda_n^3$ with 1's on the main diagonal. Note that
\begin{equation}\label{2.7}
 ps[\Lambda_n^3]=ps[\overline{\Lambda}_n^3].
\end{equation}
Indeed, it is well known that every $\Lambda_n^3$-matrix $A$ has a diagonal of ones (i.e., a set of 1's no two in the same row or column). Let $l$ be such a diagonal. There exists a permutation of rows and columns $\pi$ such that $\pi(l)$ will be the main diagonal of $\pi(A).$ Nevertheless, $per(\pi(A))=perA$ and (\ref{2.7}) follows.\newline
\indent 3) A known explicit formula for $|\overline{\Lambda}_n^3|$ (Shevelev \cite{19}) has a close structure to (\ref{2.4}):
\begin{equation}\label{2.8}
|\overline{\Lambda}_n^3|=\sum_{k=0}^{\lfloor\frac{n}{2}\rfloor}\binom{n} {k}S_{n-k}S_kU_{n-2k},
\end{equation}
\newpage
where sequence $\{S_n\}$ is defined by recursion
\begin{equation}\label{2.9}
S_0=1,\enskip S_1=0,\enskip S_n=(n-1)(S_{n-1}+\frac{1} {2} S_{n-2}), \enskip n\geq2.
\end{equation}
\indent 4) Asymptotic formula for $|\overline{\Lambda}_n^3|$  (Shevelev \cite{19})
\begin{equation}\label{2.10}
|\overline{\Lambda}_n^3|=C\sqrt{n}(\frac{n}{e})^{2n}(1+O(n^{-1+\varepsilon}),
\end{equation}
where
$$C=2\sqrt{\pi e^{-5}}=0.29098...$$
and $\varepsilon>0$ is arbitrary small for sufficiently large $n.$ \newline
\indent 5) Denote $\widehat{\Lambda}_n^3$ the set of symmetric matrices in $\overline{\Lambda}_n^3.$  $P$-spectrum on $\widehat{\Lambda}_n^3$ is given by the following theorem   (Shevelev \cite{16})
\begin{theorem} \label{t1}  Let $R(n;3)$ denote the set of all partitions of $n$ with parts more than or equal to $3.$ To every partition $r\in R(n;3): n=n_1+n_2+...+n_m, \enskip m=m(r),$ put in a correspondence the number
\begin{equation}\label{2.11}
H(r)=\prod_{i=1}^m a(n_i),
\end{equation}
where sequence $\{a(n)\}$ is defined by the recursion
\begin{equation}\label{2.12}
a(3)=6, \enskip a(4)=9,\enskip a(n)=a(n-1)+a(n-2)-2, \enskip n\geq5.
\end{equation}
Then we have
\begin{equation}\label{2.13}
 ps[\widehat{\Lambda}_n^3]=\{H(r):r\in R(n;3)\}.
 \end{equation}
\end{theorem}
\indent 6) The maximal value $M(n)$ of permanent in $\Lambda_n^3$ was found by Merriell \cite{9}.
\begin{theorem} \label{t2} If $n\equiv h\pmod 3,\enskip h=0,1,2,$ then
\begin{equation}\label{2.14}
 M(n)=6^{\frac{n-h}{3}}\lfloor(\frac{3}{2})^h\rfloor.
 \end{equation}
\end{theorem}
Note that, the case $h=0$ of (\ref{2.14}) easily follows from a general Minc-Bregman inequality for permanent of (0,1)-matrices (see \cite{11}, point 6.2, and \cite{4}).\newline
\indent 7) Put $M(n)=M^{(1})(n).$ In case of $n\equiv 0\pmod 3,$ Bolshakov \cite{3} showed that the second maximal $M^{(2)}<M^{(1})(n)$ of permanent in $\Lambda_n^3$ (such that interval $(M^{(2)},\enskip M^{(1)})$ is free from values of permanent in $\Lambda_n^3$) equals to
\begin{equation}\label{2.15}
 M^{(2)}(n)=\begin{cases}20,& if \;\; n=6,
\\120,& if\;\; n=9\\\frac{9}{16}6^{\frac{n}{3}},& if\;\;n\geq12.\end{cases}
 \end{equation}
Note that both $M^{(1)}(n)$ and $M^{(2)}(n)$ are attained in  $\widehat{\Lambda}_n^3$ (Shevelev \cite{16}).\newpage
\indent 8) Denote $m(n)$ the minimal value of permanent in $\Lambda_n^3.$ In 1979, Voorhoeve \cite{26} obtain a beautiful lower estimate for $m(n):$
\begin{equation}\label{2.16}
m(n)\geq6(\frac{4}{3})^{n-3}.
 \end{equation}
This estimate remains the best even after proof by Egorychev \cite{7} and Falikman \cite{8} the famous Van der Waerden conjectural lower estimate $per A\geq\frac {n!}{n^n}$ for every $n\times n$ doubly stochastic matrix $A.$ Indeed, this estimate yields only $m(n)\geq 3^n\frac {n!}{n^n},$ such that (\ref{2.16}) is stronger for $n\geq4.$\newline
\indent 9) Bolshakov \cite{2} found $p$-spectrum in $\Lambda_n^3$ in cases $n\leq8.$ Namely, he added to the evident $p$-spectrums
$ ps[\Lambda_3^3]=\{6\}$ and $ ps[\Lambda_4^3]=\{9\}$ also the following $p$-spectrums
 $$ ps[\Lambda_5^3]=\{12,\enskip 13\}, \enskip  ps[\Lambda_6^3]=\{17,\enskip 18,\enskip 20, \enskip 36\},$$
 \begin{equation}\label{2.17}
 ps[\Lambda_7^3]=\{24,\enskip 25,\enskip 26, \enskip 27,\enskip30,\enskip31,\enskip32,\enskip54 \},
\end{equation}
 $$ps[\Lambda_8^3]=\{33, 34, 35, 36, 37, 38, 39, 40, 41, 42, 44, 45, 48, 49, 52, 72, 78, 81 \}.$$

\section{A generalization of Theorem 1 on  matrices of class $\Lambda_n(\alpha,\beta,\gamma)$ with symmetric positions of elements}
Denote $\overline{\Lambda}_n(\alpha,\beta,\gamma)$ the set of matrices in  $\Lambda_n(\alpha,\beta,\gamma)$ with $\beta $'s on the main diagonal. It is clear that, together with (\ref{2.3}),
\begin{equation}\label{3.1}
|\overline{\Lambda}_n(\alpha,\beta,\gamma)|=K_n.
\end{equation}
Note that, as for
sets $\Lambda_n^3,\enskip \overline{\Lambda}_n^3,$ we have
\begin{equation}\label{3.2}
 ps[\Lambda_n(\alpha,\beta,\gamma)]=ps[\overline{\Lambda}_n(\alpha,\beta,\gamma)].
\end{equation}
 Denote, furthermore, $\widehat{\Lambda}_n(\alpha,\beta,\gamma)$ the set of  matrices $M=\{m_{i,j}\}$ in $\overline{\Lambda}_n(\alpha,\beta,\gamma)$ with symmetric positions of elements: $m_{i,j}=\alpha$ if and only if
 $m_{j,i}=\gamma.$\newline
 \indent $P$-spectrum on $\widehat{\Lambda}_n(\alpha,\beta,\gamma)$ is given by the following theorem.
\begin{theorem}\label{t3} If to every partition $r\in R(n;3): n=n_1+n_2+...+n_m, \enskip m=m(r),$ corresponds the number
\begin{equation}\label{3.3}
H_{\alpha,\beta,\gamma}(r)=\prod_{i=1}^m a(n_i),
\end{equation}
where sequence $\{a(n)=a(\alpha,\beta,\gamma;\enskip n\}$ is defined by the recursion
$$a(3)=\alpha^3+\beta^3+\gamma^3+3\alpha\beta\gamma,$$
$$a(4)=\alpha^4+\beta^4+\gamma^4+4\alpha\beta^2\gamma+2(\alpha\gamma)^2,$$
$$ a(n)=\beta a(n-1)+\alpha\gamma a(n-2)+$$
\begin{equation}\label{3.4}
\alpha^{n-1}(\alpha-\beta-\gamma)+\gamma^{n-1}(\gamma-\beta-\alpha), \enskip n\geq5,
\end{equation}
\newpage
then we have
\begin{equation}\label{3.5}
 ps[\widehat{\Lambda}_n(\alpha,\beta,\gamma)]=\{H_{\alpha,\beta,\gamma}(r):r\in R(n;3)\}.
 \end{equation}
\end{theorem}
\bfseries Proof. \mdseries Let $S_n$ be the symmetric permutation group of elements $1,...,n.$ Two positions $(i_1,j_1), (i_2,j_2)$ are called \slshape independent \upshape \enskip if $i_k\neq j_k,$ k=1,2. We shall say that in the $n\times n$  matrix $M=\{m_{ij}\}$ a weight $m_{ij}$ is appropriated to the position $(i,j).$ Let  $s\in S_n$ has not any cycle of length less than $n.$ Consider a map
$$\sigma: (i,j)\mapsto(s^{i}(1),\enskip s^j(1)),$$
appropriating to the position $(s^i(1), \enskip s^j(1)$ the weight  $m_{ij}.$
\begin{lemma}\label{L1}
$1)$ the map $\sigma$ is bijective;
$2)$ if $E$ is a set of pairwise independent positions, then $\sigma(E)$ is also a set of pairwise independent positions.
\end{lemma}
\bfseries Proof. \mdseries a) Consider two distinct positions
 \begin{equation}\label{3.6}
 (i_1, j_1),\enskip (i_2, j_2),
 \end{equation}
 such that, at least, one of two inequalities holds
 \begin{equation}\label{3.7}
 i_1\neq i_2,\enskip j_1\neq j_2.
 \end{equation}
 Let $ i_1\neq i_2$ such that, say, $ i_1>i_2.$ Show that $s^{i_1}(1)\neq s^{i_2}(1).$ Indeed, if to suppose that
$s^{i_1}(1)= s^{i_2}(1),$ then $s^{i_1-i_2}=1,$ i.e., $s$  has a cycle of length $i_1-i_2<n$ in spite of the condition. Conversely, if $s^{i_1}(1)\neq s^{i_2}(1),$ then $ i_1\neq i_2,$ since $s^{-1}$ has not any cycle of length less than $n$ as well.\newline
b) Let positions (\ref{3.6}) be independent. The both of inequalities (\ref{3.7}) hold and, as in a), we have $s^{i_1}(1)\neq s^{i_2}(1),\enskip s^{j_1}(1)\neq s^{j_2}(1),$ i.e. the positions $\sigma((i_1, j_1)), \sigma((i_1, j_1))$  are independent as well.\newline $\blacksquare$
\begin{lemma}\label{L2}
Let $s\in S_n$ have not any cycle of length less than $n.$ Then $(0,1)$-matrix $S$ having 1's on only positions
$$ (s^1(1),s^2(1)), (s^2(1), s^3(1)),...,(s^{n-1}(1),s^n(1)), (s^n(1),s^1(1))$$
is a incidence matrix of $s.$
\end{lemma}
\bfseries Proof. \mdseries Since $s$ has not cycles of length less than $n,$ then $\{s^1(1),...,s^n(1)\}$ is a permutation of numbers $\{1,...,n\}.$ Thus the set of positions of 1's of matrix $S$ coincides with the set of 1's of the incidence matrix of $s:\enskip (1,s(1)),...,(n,s(n)). \newline\blacksquare $ \newline
Let $P=P_n$ be $(0,1)$-matrix with 1's on positions $(1,2), (2,3),...,(n-1,n), (n,1)$ only.
\newpage
\begin{lemma}\label{L3}
 Let $s\in S_n$ have not any cycle of length less than $n.$ If $S$ and $S^{-1}$ are the incidence matrices of $s$ and $s^{-1},$ then we have
\begin{equation}\label{3.8}
\sigma^{-1}(S)=P, \enskip \sigma^{-1}(S^{-1})=P^{-1}.
 \end{equation}
\end{lemma}
\bfseries Proof. \mdseries Both of formulas follows from Lemma \ref{L2}.\newline$\blacksquare $\newline
Noting that $\sigma(I)=I$, where $I$ is the identity matrix, we conclude that
\begin{equation}\label{3.9}
S^{-1}+I+S=\sigma(P^{-1}+I+P).
 \end{equation}
Moreover, since, by the bijection $\sigma,$ to every diagonal (i.e., to every set of $n$ pairwise independent positions) of the matrix $\alpha S^{-1}+\beta I+\gamma S$ corresponds one and only one diagonal
 of the matrix $\alpha P^{-1}+\beta I+\gamma P$ with the same products of weights, then we have
 \begin{equation}\label{3.10}
per(\alpha S^{-1}+\beta I+\gamma S)=per(\alpha P^{-1}+\beta I+\gamma P).
 \end{equation}
 Note that from the definition it follows that, for every matrix $M\in \widehat{\Lambda}_n(\alpha,\beta,\gamma),$ we have a representation
  \begin{equation}\label{3.11}
M=\alpha S^{-1}+\beta I_n+\gamma S,
 \end{equation}
 where $S$ is the incidence matrix of a substitution $s.$ In case when $s$ has not any cycle of length less than n,
 the matrix $M$ is completely indecomposable matrix in $\widehat{\Lambda}_n(\alpha,\beta,\gamma).$ Thus, by (\ref{3.10}), all completely indecomposable matrices of $\widehat{\Lambda}_n(\alpha,\beta,\gamma)$ have the same permanent, equals to $per(\alpha P^{-1}+\beta I_n+\gamma P).$ \newline
 \indent In general, a substitution $s$ with the incidence matrix $S$ in (\ref{3.11}) cannot have cycles of length less than 3. Indeed, if for some $i,$ we have either $s(i)=i$ or $s(s(i))=i,$ then in both cases $s(i)=s^{-1}(i)$ which means coincidence of positions 1's of the matrices $S$ and $S^{-1}$ in the $i$-th row. \newline
 \indent Let $s\in S_n$ be an arbitrary substitution with cycles of length more than 2. Let
  \begin{equation}\label{3.12}
s=\prod_{j=1}^r s_j,
 \end{equation}
where $s_j\in S_{l_j},\enskip l_j\geq3, \sum_{j=1}^rl_j=n,$ be the decomposition of $s$ in a product of cycles. Then the matrix $M=\alpha S^{-1}+\beta I_n+\gamma S$ is a direct sum of the matrices $M_j=\alpha S^{-1}_{l_j}+\beta I+\gamma S_{l_j}$ such that, by (\ref{3.10}), $per M_j=per(\alpha P^{-1}+\beta I_{l_j}+\gamma P)$ and we have
\begin{equation}\label{3.13}
per M=\prod_{j=1}^r per M_j=\prod_{j=1}^r per (\alpha P^{-1}+\beta I_{l_j}+\gamma P).
\end{equation}
It is left to notice that Minc \cite{10} found a recursion (\ref{3.4}) for $per (\alpha I_n+\beta P+\gamma P^2)$ and, as well known, the multiplication an $n\times n$ matrix by $P^{-1}$ does not change its permanent.
\newpage
Therefore, $per (\alpha P^{-1}+\beta I_{l_j}+\gamma P)=per (\alpha I_n+\beta P+\gamma P^2).$ \newline $\blacksquare$
\begin{example}\label{e1}
Let us find $ ps[\widehat{\Lambda}_{11}(-1,3,2)].$
\end{example}
We have the following partitions of 11 with the parts not less than 3:
$$11=8+3=7+4=6+5=3+4+4=3+3+5. $$
According to (\ref{3.4}), for $a(n)=a(-1,3,2;\enskip n),$ we have
$a(3)=16,\enskip a(4)=34$ and for $n\geq5,$
$$a(n)=3a(n-1)-2a(n-2)+6(-1)^n.$$
Using induction, we find
$$ a(n)=\begin{cases} 2^{n+1},& if\enskip n \enskip is\enskip odd, \\2^{n+1}+2, & if\enskip n\enskip is\enskip even.\end{cases}$$
Therefore,
$$ps[\widehat{\Lambda}_{11}(-1,3,2)]=$$\newline$$\{a(11),\enskip a(3)a(8),\enskip a(5)a(6),\enskip a^2(3)a(5),\enskip a(3)a^2(4)\}=  $$$$ \{4096,\enskip 8224,\enskip 8320, \enskip 8704,\enskip 16384, \enskip 18496\}.$$
 $\blacksquare$\newline
In the following examples we calculate $p$-spectrum for arbitrary $n.$
\begin{example}\label{e2}
Let us find $ ps[\widehat{\Lambda}_{n}(-1,2,1)].$
\end{example}
By induction, for $a(n)=a(-1,2,1;\enskip n),$ we have
$$ a(n)=\begin{cases} 2,& if\enskip n \enskip is\enskip odd, \\4, & if\enskip n\enskip is\enskip even.\end{cases}$$
Further, again using induction, one can find that, if $n$ is even, then
$$ps[\widehat{\Lambda}_{n}(-1,2,1)]=\{4,\enskip 4^2,...,4^{\lfloor(\frac{n}{4})\rfloor} \} $$
and, if $n$ is odd, then
$$ps[\widehat{\Lambda}_{n}(-1,2,1)]=\{2,\enskip 2\cdot4,\enskip 2\cdot4^2,...,2\cdot4^{\lfloor(\frac{n-3}{4})\rfloor} \}. $$
 $\blacksquare$\newline
\begin{example}\label{e3}
Analogously, in case of $\widehat{\Lambda}_{n}(-1,1,1),$ for $a(n)=a(-1,1,1;\enskip n),$ we have
\end{example}
\newpage
$$ a(n)=\begin{cases} 4,\enskip if\enskip n\equiv0\pmod6,\\-2, \enskip if\enskip n\equiv3\pmod6,\\1,\enskip otherwise. \end{cases}$$ and
$$ps[\widehat{\Lambda}_{n}(-1,1,1)]=$$ $$\begin{cases}\{1,\enskip -2,\enskip 4,...,(-2)^{\lfloor(\frac{n-3}{3})\rfloor} \},\enskip if\enskip n\equiv1,2\pmod3,\\\{1,\enskip -2,\enskip 4,...,(-2)^{\lfloor(\frac{n-6}{3})\rfloor} ,\enskip (-2)^{\lfloor(\frac{n}{3})\rfloor}\},\enskip if\enskip n\equiv0\pmod3.\end{cases} $$
It is interesting that, in case of $n$ multiple of 3, the permanent omits the value $(-2)^{\lfloor(\frac{n-3}{3})\rfloor}.$\newline
 $\blacksquare$
 \section{Merriell type theorems  in a subclasses of $\widehat{\Lambda}_n(\alpha,\beta,\gamma)$}
Note that in class $\Lambda_n(\alpha,\beta,\gamma)$ the Minc-Bregman inequality and the Merriell theorem , generally speaking, do not hold even for positive $\alpha,\beta,\gamma.$ Nevertheles, some restrictions on $\alpha,\beta,\gamma$ allow to prove some analogs of the Merriell theorem. Recall that $M(n)$ (\ref{2.14}) is attained in $\widehat{\Lambda}_n^3.$ Denote $M_n(\alpha,\beta,\gamma)$ the maximal value of permanent in $\widehat{\Lambda}_n(\alpha,\beta,\gamma).$
\begin{theorem} \label{t4}
Consider a class  $\widehat{\Lambda}_n(\alpha,\beta,\gamma)$ with the numbers $\alpha,\beta,\gamma$ satisfying  "triangle inequlities"
\begin{equation}\label{4.1}
0\leq \alpha\leq \beta+\gamma, \enskip 0\leq \gamma\leq \alpha+\beta,
 \end{equation}
and the following additional conditions
\begin{equation}\label{4.2}
a^3(4)\leq a^4(3), \enskip \alpha\gamma+\beta (a(3))^{\frac{1}{3}}\leq (a(3))^{\frac{2}{3}},
 \end{equation}
 where sequence $\{a(n)\}$ is defined by recursion $(\ref{3.4}).$ Then, for $n$ multiple of $3,$ we have
\begin{equation}\label{4.3}
M_n(\alpha,\beta,\gamma)=(a(3))^{\frac{n}{3}}.
\end{equation}
\end{theorem}
\bfseries Proof. \mdseries Note that conditions (\ref{4.1})-(\ref{4.2}) are satisfied, e.g., in case $\alpha=\beta=\gamma=1.$ Using induction, let us prove that
\begin{equation}\label{4.4}
a(n)\leq (a(3))^{\frac{n}{3}}.
\end{equation}
Indeed, for $n=3,$ this inequality is trivial, while, for $n=4,$ it follows from the first condition (\ref{4.2}). Let it hold for $n\leq m-1.$ Then, according to (\ref{3.4}), we have
$$a(m)=\beta a(m-1)+\alpha\gamma a(m-2)+\alpha^{m-1}(\alpha-\beta-\gamma)+\gamma^{m-1}(\gamma-\alpha-\beta)\leq$$
$$\beta(a(3))^{\frac{m-1}{3}}+\alpha\gamma(a(3))^{\frac{m-2}{3}}=$$ $$(a(3))^{\frac{m-2}{3}}(\beta(a(3))^{\frac{1}{3}}+\alpha\gamma)\leq(a(3))^{\frac{m-2}{3}}(a(3))^{\frac{2}{3}}=
(a(3))^{\frac{m}{3}}.$$
\newpage
Note that, according Theorem \ref{t3}, the equality in (\ref{4.4}) holds in a direct sum of $(3\times3)$-matrices of $\widehat{\Lambda}_3(\alpha,\beta,\gamma)$ which corresponds to the partition $n=3+3+...+3.$
\indent Let now $A\in \widehat{\Lambda}_n(\alpha,\beta,\gamma).$ By Theorem \ref{t3}, there exists a partition of $n$ with the parts not less than 3, $n=n_1+...+n_m,$ such that
$$ per A=\prod_{i=3}^m a(n_i)$$
and, in view of (\ref{4.4}), we have
$$ per A\leq\prod_{i=3}^m a(3)^{\frac{n_i}{3}}=(a(3))^{\frac{n}{3}}.$$
This proves (\ref{4.3}).\newline
$\blacksquare$
\begin{example}\label{e4}
Consider case $\beta=\gamma-\alpha.$
\end{example}
Let us find the values of $\alpha,$ depending on the magnitude of $\gamma,$ for which the conditions of Theorem \ref{t4} are satisfied. According to $(\ref{3.4}),$ we have
\begin{equation}\label{4.5}
 a(3)=\alpha^3+(\gamma-\alpha)^3+\gamma^3+3\alpha(\gamma-\alpha)\gamma=2\gamma^3,
\end{equation}
\begin{equation}\label{4.6}
 a(4)=\alpha^4+(\gamma-\alpha)^4+\gamma^4+4\alpha(\gamma-\alpha)^2\gamma + 2(\alpha\gamma)^2=2(\alpha^4+\gamma^4).
\end{equation}
Thus the condition $a^3(4)\leq a^4(3)$ means that $\alpha^4+\gamma^4\leq2^{\frac{1}{3}}\gamma^4,$ or
\begin{equation}\label{4.7}
0\leq\alpha\leq (2^{\frac{1}{3}}-1)^{\frac{1}{4}}\gamma=0.7140199...\gamma.
\end{equation}
and it is easy to verify that the second condition in $(\ref{4.2})$ is satisfied as well. As a collorary, we obtain the following result.
\begin{theorem}\label{t5}
If $(\ref{4.7})$ holds, then, for $n$ multiple of $3,$ we have
\begin{equation}\label{4.8}
M_n(\alpha,\gamma-\alpha,\gamma)=2^{\frac{n}{3}}\gamma^n.
\end{equation}
\end{theorem}
$\blacksquare$\newline
\indent Simple forms of sequence $\{a(n)\}$ in Examples \ref{e1}-\ref{e2} allow to suppose that in case $\beta=\gamma-\alpha$ (or symmetrical case $\beta=\alpha-\gamma$) sequence $\{a(n)\}$ keeps a sufficiently simple form. We find this form in the following lemma.
\begin{lemma}\label{L4}
If $\beta=\gamma-\alpha,$ then sequence $\{a(n)\}$ which is defined by recursion $(\ref{3.4})$ has the form
\begin{equation}\label{4.9}
 a(n)=\begin{cases} 2\gamma^n,& if\enskip n \enskip is\enskip odd, \\2(\alpha^n+\gamma^n), & if\enskip n\enskip is\enskip even.\end{cases}
\end{equation}
\end{lemma}
\bfseries Proof. \mdseries
Using induction with the base $(\ref{4.5})$ -$(\ref{4.6}),$ suppose that $(\ref{4.9})$ holds for $m\leq n.$ Then, by  $(\ref{3.4}),$ for even $n,$ we have
$$a(n+1)=(\gamma-\alpha)a(n)+\alpha\gamma a(n-1)+2\alpha^{n}(\alpha-\gamma)=$$
\newpage
$$ 2(\gamma-\alpha)(\alpha^n+\gamma^n)+2\alpha\gamma^{n}+2\alpha^{n}(\alpha-\gamma)=2\gamma^{n+1}, $$

while, if $n$ is odd, then we have
$$a(n+1)=2(\gamma-\alpha)\gamma^n+$$ $$2\alpha\gamma(\alpha^{n-1}+\gamma^{n-1})+2\alpha^n(\alpha-\gamma)
=2(\alpha^{n+1}+\gamma^{n+1}). $$
$\blacksquare$\newline
\indent Let now
\begin{equation}\label{4.10}
\alpha\geq (2^{\frac{1}{3}}-1)^{\frac{1}{4}}\gamma=0.7140199...\gamma.
\end{equation}
\begin{theorem}\label{t6}
If $(\ref{4.10})$ holds, then, for $n$ multiple of $4,$ we have
\begin{equation}\label{4.11}
M_n(\alpha,\gamma-\alpha,\gamma)=(2(\alpha^4+\gamma^4))^{\frac{n}{4}}.
\end{equation}
\end{theorem}
\bfseries Proof. \mdseries  From  $(\ref{4.5}),  (\ref{4.6})$ and (\ref{4.10})  we conclude that
\begin{equation}\label{4.12}
a(3)\leq (a(4))^{\frac{3}{4}}.
\end{equation}
Let us show that, for $n\geq3,$
\begin{equation}\label{4.13}
a(n)\leq (a(4))^{\frac{n}{4}}.
\end{equation}
For $n=4,$ inequality (\ref{4.13}) is trivial. For $n\geq5,$ we have
$$\alpha^n+\gamma^n=\alpha^n(1+(\frac{\gamma}{\alpha})^n)\leq \alpha^n(1+(\frac{\gamma}{\alpha})^4)^{\frac{n}{4}}\leq(\alpha^4+\gamma^4)^{\frac{n}{4}}$$
and thus, using Lemma \ref{L4}, we have
$$a(n)\leq2(\alpha^n+\gamma^n)<2^{\frac{n}{4}} (\alpha^4+\gamma^4)^{\frac{n}{4}}=(a(4))^{\frac{n}{4}}, \enskip n\geq3.$$
Let now $A\in \widehat{\Lambda}_n(\alpha,\beta,\gamma).$ By Theorem \ref{t3}, there exists a partition of $n$ with the parts not less than 3, $n=n_1+...+n_m,$ such that
$$ per A=\prod_{i=3}^m a(n_i)$$
and, in view of (\ref{4.13}), we have
$$ per A\leq\prod_{i=3}^m a(4)^{\frac{n_i}{4}}=(a(3))^{\frac{n}{4}}$$
with the equality in a direct sum of $(4\times4)$-matrices of $\widehat{\Lambda}_3(\alpha,\beta,\gamma)$ which corresponds to the partition $n=4+4+...+4.$\newline
$\blacksquare$\newline
Note that, if $\alpha\neq (2^{\frac{1}{3}}-1)^{\frac{1}{4}}\gamma,$ then in Theorem 5 we have only maximizing
matrix (up to a permutation of the rows and columns) which corresponds to the partition $n=3+3+...+3;$ in Theorem \ref{t6}
we also have only maximizing matrix (up to a permutation of the rows and columns) which corresponds to the partition $n=4+4+...+4.$ It is interesting that, only in case of the equality \newpage $\alpha=\theta\gamma,$ where $\theta=(2^{\frac{1}{3}}-1)^{\frac{1}{4}},$ the both of Theorems 5-\ref{t6} are true for every $n$ multiple of 12 with the equality of the maximums: $(2(\alpha^4+\gamma^4))^{\frac{n}{4}}=2^{\frac{n}{3}}\gamma^n.$ Thus, up to a positive factor $\gamma,$ \slshape the class
\begin{equation}\label{4.14}
\widehat{\Lambda}_n(\theta,1-\theta,1), \enskip \theta=(2^{\frac{1}{3}}-1)^{\frac{1}{4}},
\end{equation}
possesses an interesting extremal property: it contains $\frac {n}{12}+1$ maximizing matrices (up to a permutation of the rows and columns), instead of only maximizing matrix, if $\theta\neq(2^{\frac{1}{3}}-1)^{\frac{1}{4}}.$\upshape \newline
\indent Indeed, the number of the maximizing matrices (up to a permutation rows and columns), is defined by the number of the following partitions of $n\equiv\pmod {12}:$
$$n=3+3+...+3,\enskip n=(4+4+4)+(3+...+3),...,$$
 $$ n=\underbrace{4+...+4}_{3i}+\underbrace{3+...+3}_{n-12i},\enskip i=0,1,...,\frac{n}{12}.$$
 $\blacksquare$
\section{Estimate of cardinality of p-spectrum on circulants in $\Lambda_n^3$ and ${\Lambda}_n(\alpha,\beta,\gamma)$}
Denote $\Delta_n^3\subset \Lambda_n^3$ the set of the circulants in $\Lambda_n^3.$ Note that a circulant $A\in \Delta_n^3$ has a form $A=P^i+P^j+P^k,\enskip 0\leq i<j<k\leq n,$ where $P=P_n$ is $(0,1)$-matrix with 1's on positions $(1,2), (2,3),...,(n-1,n), (n,1)$ only. Multiplicating $A$ by $P^{-i}$ we obtain circulant $B$ of the form  $$B=I_n+P^r+P^s$$
with $per B=per A.$ Since $B$ is defined by a choice of two different values  $0<r<k\leq n,$ then trivially $$ps[\Delta_n^3]\leq \binom {n} {2}<\frac{n^2}{2}.$$
Now we prove essentially more exact and practically unimprovable estimate.
\begin{theorem}\label{t7}
We have
\begin{equation}\label{5.1}
ps[\Delta_n^3]\leq\lfloor\frac{n^2+3}{12} \rfloor.
\end{equation}
\end{theorem}
\bfseries Proof. \mdseries Let us back to the general form
$$A=P^i+P^j+P^k,\enskip 0\leq i<j<k\leq n.$$
Note that, $A$ is defined by a choice of a vector $(i,j,k),$ but its rotation, i.e., a passage to a vector of the form
$(i+l, j+l,k+l) \pmod {n},$ does not change the magnitude of $per A.$ Indeed it corresponds to the multiplication $A$ by $P^l,$ and our statement follows from the equality $per (P^lA)=per A.$ Besides,\newpage its reflection relatively some diameter of the imaginary circumference of the rotation, by the symmetry, keeps magnitude of the permanent. Since geometrically three points on the imaginary circumference define a triangle, then our problem reduces to a triangle case of the following general problem, posed by Professor Richard H. Reis (South-East University of Massachusetts, USA) in a private communication to Hansraj Gupta in 1978): \newline \indent" Let a circumference is split by the same n parts. It is required to find the number $R(n,k)$ of the incongruent convex $k$-gons, which could be obtaind by connection of some $k$ from $n$ dividing points. Two $k$-gons are considered congruent if they are coincided at the rotation of one relatively other along the circumference and (or) by reflection of one of the $k$-gons relatively some diameter."\newline
\indent In 1979, Gupta \cite{6} gave a solution of the Reis problem in the form (for a short solution, see author's paper \cite{22}):
\begin{equation}\label{5.2}
R(n,k)=\frac{1}{2}(\binom{\lfloor(\frac {n-h_k}{2})\rfloor} {\lfloor(\frac{k} {2})\rfloor}+\frac{1}{k}\sum_{d|\gcd(k,n)}\varphi(d)\binom{\frac{n}{d}-1}{\frac{k}{d}-1}).
\end{equation}
If to denote $\Delta_n^k\subset \Lambda_n^k$ the set of the circulants in $\Lambda_n^k,$ then from our arguments it follows that
\begin{equation}\label{5.3}
ps[\Delta_n^k]\leq R(n,k).
\end{equation}
In case $k=3,$ from (\ref{5.2})-(\ref{5.3}) we find
$$ps[\Delta_n^3]\leq \begin{cases} \frac{n^2}{12},\enskip if\enskip n\equiv0\pmod6,\\ \frac{n^2-1}{12}, \enskip if\enskip n\equiv1,5\pmod6,\\ \frac{n^2-4}{12},\enskip if\enskip n\equiv2,4\pmod6,\\ \frac{n^2+3}{12},\enskip if\enskip n\equiv3\pmod6, \end{cases}$$
and (\ref{5.1}) follows.\newline
$\blacksquare$
\begin{example}\label{e5}
In case $n=5$ we have only two incongruent triangles corresponding to circulants $I_5+P+P^2$ and $I_5+P+P^3.$
\end{example}
Nevertheless, the calculations give $per(I_5+P+P^2)=per(I_5+P+P^3)=13.$ Thus $ps[\Delta_5^3]=\{13\},$ and $|ps[\Delta_5^3]|=1.$
\begin{example}\label{e6}
In case $n=6$ we have three incongruent triangles corresponding to circulants $I_6+P+P^2,\enskip I_6+P+P^3$ and $I_6+P^2+P^4.$
 \end{example}
 The calculations give $per(I_6+P+P^2)=20,\enskip per(I_6+P+P^3)=17,$ while $per(I_6+P^2+P^4)=36.$ Thus $ps[\Delta_6^3]=\{17,20,36\},$ and $|ps[\Delta_6^3]|=3.$\newline
 \indent Note that a respectively large magnitude of $per(I_6+P^2+P^4)$ is explained\newpage by its decomposability in a direct product of circulants  $(I_3+P+P^2)\otimes(I_3+P+P^2),$ such that $per(I_6+P^2+P^4)=(per (I_3+P+P^2))^2=6^2=36.$\newline
\indent It is clear that, in case of circulants in ${\Lambda}_n(\alpha,\beta,\gamma),$ the upper estimate (\ref{5.1}) yields either the same estimate, if $\alpha=\beta=\gamma,$ or $\lfloor\frac{n^2+3}{4} \rfloor, $ if $\alpha=\beta\neq\gamma$ (and in the symmetric cases), or  $\lfloor\frac{n^2+3}{2} \rfloor, $ if $\alpha,\beta,\gamma$ are distinct numbers.\newline
\indent Add that a bijection indicated in \cite{22} allows to apply formula (\ref{5.2}) to
enumerating the two-color bracelets of $n$ beads, $k$ of which are black and $n-k$ are white
(see, e.g., the author's explicit formulas for sequences A032279-A032282,\enskip A005513-A005516 in \cite{23}).

\section{Algorithm of calculations of upper magnitudes of the permanent in $\widehat{\Lambda}_n^3$}
Theorem \ref{t1} allows, using some additional arguments, to give an algorithm of calculations of upper magnitudes of the permanent in $\widehat{\Lambda}_n^3.$ In connection with this, we need an important lemma for numbers (\ref{2.12}).
\begin{lemma} \label{L5}
For $n_1, n_2\geq3,$ we have
\begin{equation}\label{6.1}
a(n_1+n_2)\leq a(n_1)a(n_2)\enskip , n_1, n_2\geq3.
 \end{equation}
\end{lemma}
\bfseries Proof. \mdseries
 By usual way, from (\ref{2.12}) we find
\begin{equation}\label{6.2}
a(n)=\varphi^n+2+(-1)^n\varphi^{-n},
\end{equation}
where $\varphi=\frac{\sqrt{5}+1}{2}$ is the golden ratio.\newline
Denote $\varepsilon(n)=(-1)^n\varphi^{-n}.$  Since $n\geq3,$ then $|\varepsilon(n)|<0.24,$ and, consequently, if $n=n_1+n_2,$ then $(2+\varepsilon(n_1))(2+\varepsilon(n_2)>1.76^2>3.$ Therefore, we have
$$a(n_1)a(n_2)=(\varphi^{n_1}+2+\varepsilon(n_1))(\varphi^{n_2}+2+\varepsilon(n_2))> $$
\begin{equation}\label{6.3}
 \varphi^{n_1+n_2}+3>\varphi^{n_1+n_2}+2+\varepsilon(n_1+n_2)=a(n_1+n_2).
\end{equation}
$\blacksquare$\newline
\indent Note that, actually, the difference between the hand sides of ({6.3}) more than $1.76(\varphi^{n_1}+\varphi^{n_2}).$ \newline
\indent Let now $n\equiv j\pmod{3}, \enskip j=0,1,2,$ and $t\in \mathbf{N}.$ Let $R(m;\nu)$ denote the set of all partitions of $n$ with parts more than or equal to $\nu.$ For us an important role play cases $\nu=3,4.$ To $r\in R(m;3),\enskip \rho\in R(m;4)$ put in a correspondence the sets
\begin{equation}\label{6.4}
H_{m;3}(r)=\{\Pi_{r_i\in r}a_{r_i}\}; \enskip H_{m;4}(\rho)=\{\Pi_{\rho_i\in \rho}a_{\rho_{i}}\}.
\end{equation}
In case $m=3,$ when $\rho=\varnothing,$ let us agree that $H_{3;4}$ is a singleton $\{6\}$.
 \newline \indent Consider now the set $L_t^{(j)}=L_t^{(j)}(n):$\newpage
\begin{equation}\label{6.5}
L_t^{(j)}=\bigcup_{i=1}^{4t+j}\{6^{\frac {n-j-3i}{3}}y: y\in H_{3i+j;4}(\rho), y\geq9^{3t+j}6^{i-4t-j}\}.
\end{equation}
\begin{theorem} \label{t8} (algorithm of calculation of upper magnitudes of the permanent in $\widehat{\Lambda}_n^3).$ If $n\geq4(3t+j),$ then the ordered over decrease set $L_t^{(j)}$ gives the $|L_t^{(j)}|$ upper magnitudes of the permanent in $\widehat{\Lambda}_n^3.$
\end{theorem}
\bfseries Proof. \mdseries Note that the proof is the same for every value of $j.$ Therefore, let us consider, say, $j=0.$ From ({6.1}) it follows that, if $r\in R(n;3)$ contains $\lambda_3$ parts 3 and $\lambda_3\leq \frac{n}{3}-4t,$ then, for $y\in H_{n-3\lambda_3;4}(\rho),$ we have
$$ 6^{\lambda_3}y\leq6^{\frac{n}{3}-4t}9^{3t}.$$
This means that for the formation the list of all upper magnitudes of the permanent in $\widehat{\Lambda}_n^3$ in the condition $n\geq12t,$ which are bounded from below by $6^{\frac{n}{3}-4t}9^{3t},$ it is sufficient to consider only a part of the spectrum containing numbers $\{6^{\lambda_3}y\},$ where $y\in H_{n-3\lambda_3;4}(\rho)$ with the opposite condition $\lambda_3\geq \frac{n}{3}-4t.$ From the equality $3\lambda_3+...+n\lambda_n=n$ with the condition $\lambda_3\geq \frac{n}{3}-4t,$ we have
$$4\lambda_4+...+n\lambda_n\leq n-3(\frac{n}{3}-4t)=12t.$$
Since $12t$ does not depend on $n,$ there is only a finite assembly of such partition for arbitrary $n.$ This ensures a possibility of the realization of the algorithm.\newline
\indent For the considered $r\in R(n;3),$ for $\lambda_3\geq\frac{n}{3}-4t,$ we have $H_{n;3}(r)=6^{\frac{n-m}{3}},$
where $y\in H_{m;4}(\rho),$ and $m$ has the form $m=3i,\enskip 1\leq i\leq4t.$ Thus we should choose only $H_{n;3}(r)\geq 6^{\frac {n}{3}-4t}9^{3t},$ and this yields
$$y\geq9^{3t}6^{\frac {n}{3}-4t}=9^{3t}6^{i-4t}, \enskip i=1,2,...,4t.$$
$\blacksquare$\newline
\indent In order to use Theorem {8} for calculation the upper magnitudes $M^{(1})(n)>M^{(2})(n)>...$ of the permanent in $\widehat{\Lambda}_n^3,$ in case, say, $n\equiv0\pmod3,$ \newline
 1) we write a list of partition of numbers $3i,\enskip i=2,3,...,4t$ with the parts not less than 4. \newline
 2) The corresponding values of $y$ we compare with $9^{3t}6^{i-4t}$ and keep only $y\geq 9^{3t}6^{i-4t}.$ \newline
 3) After that we regulate over decrease numbers $\{y6^{\frac{n}{3}-i}\}.$\newline
\indent Below we give the first 10 upper magnitudes $\widehat{M}^{(1)}>\widehat{M}^{(2)}>...>\widehat{M}^{(10)},$ of the permanent in $\widehat{\Lambda}_n^3$ for $n\geq24,$ via numbers $\{a(n)\}$ (\ref{2.12}).
\begin{equation}\label{6.6}
 \widehat{M}^{(1)}=\begin{cases} a(3)^{\frac{n}{3}}=6^{\frac{n}{3}},\enskip if\enskip n\equiv0\pmod3,\\a(4)a(3)^{\frac{n-4}{3}}=\frac{3}{2}6^{\frac{n-1}{3}}, \enskip if\enskip n\equiv1\pmod3,\\a(4)^2a(3)^{\frac{n-8}{3}}=\frac{9}{4}6^{\frac{n-2}{3}},\enskip if\enskip n\equiv2\pmod3. \end{cases}
\end{equation}\newpage
Formula (\ref{6.6}) shows that $M^{(1)}(n)$ is attained in  $\widehat{\Lambda}_n^3.$
\begin{equation}\label{6.7}
 \widehat{M}^{(2)}=\begin{cases}a(4)^3a(3)^{\frac{n-12}{3}}=\frac{9}{16}6^{\frac{n}{3}},\enskip if\enskip n\equiv0\pmod3,\\a(7)a(3)^{\frac{n-7}{3}}=\frac{31}{36}6^{\frac{n-1}{3}}, \enskip if\enskip n\equiv1\pmod3,\\a(5)a(3)^{\frac{n-5}{3}}=\frac{13}{6}6^{\frac{n-2}{3}},\enskip if\enskip n\equiv2\pmod3; \end{cases}
\end{equation}

\begin{equation}\label{6.8}
\widehat{M}^{(3)}=\begin{cases}a(6)a(3)^{\frac{n-6}{3}}=\frac{5}{9}6^{\frac{n}{3}},\enskip if\enskip n\equiv0\pmod3,\\a(4)^4a(3)^{\frac{n-16}{3}}=\frac{27}{32}6^{\frac{n-1}{3}}, \enskip if\enskip n\equiv1\pmod3,\\a(8)a(3)^{\frac{n-8}{3}}=\frac{49}{36}6^{\frac{n-2}{3}},\enskip if\enskip n\equiv2\pmod3; \end{cases}
\end{equation}

\begin{equation}\label{6.9}
 \widehat{M}^{(4)}=\begin{cases}a(4)a(5)a(3)^{\frac{n-9}{3}}=\frac{13}{24}6^{\frac{n}{3}},\enskip if\enskip n\equiv0\pmod3,\\a(4)a(6)a(3)^{\frac{n-10}{3}}=\frac{5}{6}6^{\frac{n-1}{3}}, \enskip if\enskip n\equiv1\pmod3,\\a(4)a(7)a(3)^{\frac{n-11}{3}}=\frac{31}{24}6^{\frac{n-2}{3}},\enskip if\enskip n\equiv2\pmod3; \end{cases}
\end{equation}

\begin{equation}\label{6.10}
 \widehat{M}^{(5)}=\begin{cases}a(9)a(3)^{\frac{n-9}{3}}=\frac{13}{36}6^{\frac{n}{3}},\enskip if\enskip n\equiv0\pmod3,\\a(4)^2a(5)a(3)^{\frac{n-13}{3}}=\frac{13}{16}6^{\frac{n-1}{3}}, \enskip if\enskip n\equiv1\pmod3,\\a(4)^5a(3)^{\frac{n-20}{3}}=\frac{81}{64}6^{\frac{n-2}{3}},\enskip if\enskip n\equiv2\pmod3; \end{cases}
\end{equation}

\begin{equation}\label{6.11}
 \widehat{M}^{(6)}=\begin{cases}a(4)a(8)a(3)^{\frac{n-12}{3}}=\frac{49}{144}6^{\frac{n}{3}},\enskip if\enskip n\equiv0\pmod3,\\a(5)^2)a(3)^{\frac{n-10}{3}}=\frac{169}{216}6^{\frac{n-1}{3}}, \enskip if\enskip n\equiv1\pmod3,\\a(4)^2a(6)a(3)^{\frac{n-14}{3}}=\frac{5}{4}6^{\frac{n-2}{3}},\enskip if\enskip n\equiv2\pmod3; \end{cases}
\end{equation}

\begin{equation}\label{6.12}
 \widehat{M}^{(7)}=\begin{cases}a(4)^2a(7)a(3)^{\frac{n-15}{3}}=\frac{31}{96}6^{\frac{n}{3}},\enskip if\enskip n\equiv0\pmod3,\\a(10)a(3)^{\frac{n-10}{3}}=\frac{125}{216}6^{\frac{n-1}{3}}, \enskip if\enskip n\equiv1\pmod3,\\a(4)^3a(5)a(3)^{\frac{n-17}{3}}=\frac{39}{32}6^{\frac{n-2}{3}},\enskip if\enskip n\equiv2\pmod3; \end{cases}
\end{equation}

\begin{equation}\label{6.13}
 \widehat{M}^{(8)}=\begin{cases}a(4)^6a(3)^{\frac{n-24}{3}}=\frac{81}{256}6^{\frac{n}{3}},\enskip if\enskip n\equiv0\pmod3,\\a(4)a(9)a(3)^{\frac{n-13}{3}}=\frac{13}{24}6^{\frac{n-1}{3}}, \enskip if\enskip n\equiv1\pmod3,\\a(5)a(6)a(3)^{\frac{n-11}{3}}=\frac{65}{54}6^{\frac{n-2}{3}},\enskip if\enskip n\equiv2\pmod3; \end{cases}
\end{equation}

\begin{equation}\label{6.14}
 \widehat{M}^{(9)}=\begin{cases}a(4)^3a(6)a(3)^{\frac{n-18}{3}}=\frac{5}{16}6^{\frac{n}{3}},\enskip if\enskip n\equiv0\pmod3,\\a(5)a(8)a(3)^{\frac{n-13}{3}}=\frac{637}{1296}6^{\frac{n-1}{3}}, \enskip if\enskip n\equiv1\pmod3,\\a(5)^2a(4)a(3)^{\frac{n-14}{3}}=\frac{169}{144}6^{\frac{n-2}{3}},\enskip if\enskip n\equiv2\pmod3; \end{cases}
\end{equation}

\begin{equation}\label{6.15}
 \widehat{M}^{(10)}=\begin{cases}a(5)a(7)a(3)^{\frac{n-12}{3}}=\frac{403}{1296}6^{\frac{n}{3}},\enskip if\enskip n\equiv0\pmod3,\\a(4)^3a(7)a(3)^{\frac{n-19}{3}}=\frac{31}{64}6^{\frac{n-1}{3}}, \enskip if\enskip n\equiv1\pmod3,\\a(11)a(3)^{\frac{n-11}{3}}=\frac{67}{72}6^{\frac{n-2}{3}},\enskip if\enskip n\equiv2\pmod3. \end{cases}
\end{equation}

\section{Main conjectural inequality for maximum of permanent in completely indecomposable $\Lambda_n^3 $-matrices}
Denote $\Lambda_{n, 1}^3$ the set of completely indecomposable matrices in  $\Lambda_n^3,$ i. e., the set of $\Lambda_n^3$-matrices containing no $\Lambda_m^3$-submatrices. Let $\mu_1 (n)$ denote the\newpage maximum of permanent in
$\Lambda_{n, 1}^3.$ Our very plausible conjecture which we call \slshape "main conjectural inequality (MCI)" \upshape \enskip is the following.
\begin{conjecture}\label{c1} (Cf. \cite{21}, pp. 165-166)
For $n_1, n_2\geq3,$ we have
\begin{equation}\label{7.1}
\mu_1 (n_1+n_2)\leq\mu_1 (n_1)\mu_1 (n_2).
 \end{equation}
\end{conjecture}
In Lemma \ref{L5} we essentially proved that in subclass $\widehat{\Lambda}_n^3$ the MCI is valid.\newline
\indent Besides, in all known cases MCI holds. Moreover, as we shall see, our algorithm of calculation the consecutive upper magnitudes $(M=M_1>M_2>...)$ of permanent in $\Lambda_n^3,$ which is based on MCI, reproduces all Merriell's and Bolshakov's results for $M_1$ and $M_2.$ \enskip Note also that, for sufficiently large $n,$ the number of consecutive upper magnitudes of permanent in $\Lambda_n^3$ grows very quickly with every step of extension of the list of known $p$-specrums for small $n.$
E.g., using the found by Bolshakov $ ps[\Lambda_i^3],\enskip i\leq8,$ we obtain, for sufficiently large $n,$ 4, 7 and 11 upper values of $ps[\Lambda_n^3]$ in cases $n=3k,\enskip 2k+1$ and $3k+2$ correspondingly. After calculation $ ps[\Lambda_9^3],$ the number of upper values of, e.g., $ps[\Lambda_{3k}^3]$ increases more than thrice.

\section{Algorithm of calculations of upper magnitudes of the permanent in $\Lambda_n^3$ based on MCI}
Let $n\equiv j\pmod{3}, \enskip j=0,1,2,$ and $t\in \mathbf{N}.$ Let $R(m;\nu)$ denote the set of all partitions of $n$ with parts more than or equal to $\nu.$ For us an important role play cases $\nu=3,4.$ To $r\in R(m;3),\enskip \rho\in R(m;4)$ put in a correspondence sets
\begin{equation}\label{8.1}
\pi_{m;3}(r)=\{\Pi_{r_i\in r}x_{r_i}\}; \enskip \pi_{m;4}(\rho)=\{\Pi_{\rho_i\in \rho}x_{\rho_{i}}\}
\end{equation}
where $x_s$ runs through all values of permanent in set $\Lambda_{s, 1}^3$ of completely indecomposable matrices in  $\Lambda_s^3$ ( in case $m=3,$ when $\rho=\varnothing,$ let us agree that $\pi_{3;4}$ is a singleton $\{6\}$). \newline \indent Consider now the set $E_t^{(j)}=E_t^{(j)}(n):$
\begin{equation}\label{8.2}
E_t^{(j)}=\bigcup_{i=1}^{4t+j}\{6^{\frac {n-j-3i}{3}}y: y\in \pi_{3i+j;4}(\rho), y\geq9^{3t+j}6^{i-4t-j}\}.
\end{equation}
\begin{theorem} \label{t9} (algorithm of calculation of upper magnitudes of the permanent in $\Lambda_n^3).$ If $n\geq4(3t+j),$ then the ordered over decrease set $E_t^{(j)}$ gives the $|E_t^{(j)}|$ upper magnitudes of the permanent in $\Lambda_n^3.$
\end{theorem}
\bfseries Proof. \mdseries We need three lemmas.
 \begin{lemma}\label{L6}
For $n\geq4,$ we have
\begin{equation}\label{8.3}
\mu_1(n)\leq3^{\frac{n}{2}}.
\end{equation}
\end{lemma}\newpage
\bfseries Proof. \mdseries Let, firstly, $n\equiv0\pmod{4}.$ Note that $\mu_1(4)=D_4=9.$ Using ({7.1}), we find
$$\mu_1(n)\leq \mu_1(4)\mu_1(n-4)\leq...\leq (\mu_1(4))^{\frac{n}{4}}=3^{\frac{n}{2}}.$$
Let, furthermore, $n\equiv i\pmod{4}, \enskip i=1,2,3.$ Note that, by (\ref{2.17}), $\mu_1(5)\leq13<3^{2.5}.$ Therefore, using ({7.1}), we have
$$\mu_1(n)\leq (\mu_1(4))^{\frac{n-5i}{4}}(\mu_1(5))^i<3^{\frac{n-5i}{2}}3^{2.5i}=3^{\frac{n}{2}}.$$
$\blacksquare$
\begin{lemma}\label{L7} Let $n=3\lambda_3+4\lambda_4+...+n\lambda_n$ be a partition of $n$ with the parts not less than $3.$ If $\lambda_3\leq l,$ and $n$ has the form $ n=3l+4m,$ then, for completely indecomposable matrices $A_i\in \Lambda_i^3, \enskip i=3,4,...,n,$ we have
\begin{equation}\label{8.4}
\prod_{i=3}^n(per A_i)^{\lambda_i}\leq 6^l9^m.
\end{equation}
\end{lemma}
\bfseries Proof. \mdseries Using Lemma \ref{L6}, we have
$$\prod_{i=3}^n(per A_i)^{\lambda_i}\leq6^{\lambda_3}\sqrt{3}^{4\lambda_4+...+n\lambda_n }\leq6^l\sqrt{3}^{n-\lambda_3}=6\sqrt{3}^{4m}=6^l9^m.$$
$\blacksquare$
\begin{lemma}\label{L8} Let $n=3\lambda_3+4\lambda_4+...+n\lambda_n$ and
$$\lambda_3\leq \frac {n-4j}{3}-4t,\enskip n\geq4(3t+j),$$
where $t$ is a nonnegative integer and $j$ is the residue of $n$ modulo $3,\enskip j=0,1,2,$ then,
for completely indecomposable matrices $A_i\in \Lambda_i^3, \enskip i=3,4,...,n,$ we have
\begin{equation}\label{8.5}
\prod_{i=3}^n(per A_i)^{\lambda_i}\leq 6^{\frac {n-4j}{3}-4t}9^{3t+j}.
\end{equation}
\end{lemma}
\bfseries Proof. \mdseries Put $l=\frac {n-4j}{3}-4t,\enskip m=\frac {n-3l}{4}=3t+j.$ Now the lemma follows from Lemma \ref{L7}.\newline
$\blacksquare$\newline
 It is left to note that, after these  lemmas, the proof of Theorem \ref{t9} is the same as proof of Theorem \ref{t8}.\newline
 $\blacksquare$\newline
 \indent Note that the using of this algorithm is based on the small elements of $p$-spectrum. \newline
 \indent Consider, e.g., case $t=0,\enskip j=2.$ According to (\ref{8.2}), we have

$$ E_0^{(2)}=\bigcup_{i=1}^{2}\{6^{\frac {n-2-3i}{3}}y: y\in \pi_{3i+j;4}(\rho), y\geq81\cdot6^{i-2}\}=$$\newpage
\begin{equation}\label{8.6}
\{6^{\frac {n-5}{3}}per A,\enskip A\in \Lambda_5^3: 6per A\geq81\}\cup\{81\cdot6^{\frac {n-8}{3}}\}.
\end{equation}
Note that, the second set is a simpleton, since, by MCI, $\mu_1(8,3)\leq(\mu_1(4,3))^2=81.$ Since, by (\ref{2.17}), $M(5)=13<\frac{81}{6},$ then the first set in (\ref{8.6}) is empty. Thus $E_0^{(2)}=E_0^{(2)}(n)$ is simpleton:
$$ E_0^{(2)}=\{81\cdot6^{\frac {n-8}{3}}\}$$

and we have
$$M^{(1)}(n)=81\cdot6^{\frac {n-8}{3}}, \enskip n\geq8,$$
which corresponds to Merriell's result in case $n\equiv2\pmod3.$\newline
Further research of the set (\ref{8.2}), using (\ref{2.17}), gives the following results:\newline
1) $j=0,\enskip n\geq24.$
$$M^{(1)}(n)=6^{\frac{n}{3}},\enskip M^{(2)}(n)=\frac{9}{16}6^{\frac{n}{3}},$$
\begin{equation}\label{8.7}
M^{(3)}(n)=\frac{5}{9}6^{\frac{n}{3}},\enskip M^{(4)}(n)=\frac{13}{24}6^{\frac{n}{3}}.
\end{equation}
The continuation of this list requires the knowing of $ps[\Lambda_9^3].$ Note that a more detailed analysis shows that after calculation $ps[\Lambda_9^3]$ in this case one can obtain the first $12+|G|$ upper magnitudes of the permanent in  $\Lambda_n^3$, where $G=ps[\Lambda_9^3]\cap ([69,116]\setminus \{72,78,102,108\}).$\newline
2) $j=1,\enskip n\geq28.$
$$M^{(1)}(n)=\frac{3}{2}6^{\frac{n-1}{3}},\enskip M^{(2)}(n)=\frac{8}{9}6^{\frac{n-1}{3}},$$

$$M^{(3)}(n)=\frac{31}{36}6^{\frac{n-1}{3}},\enskip M^{(4)}(n)=\frac{27}{32}6^{\frac{n-1}{3}},$$
\begin{equation}\label{8.8}
M^{(5)}(n)=\frac{5}{6}6^{\frac{n-1}{3}},\enskip M^{(6)}(n)=\frac{13}{15}6^{\frac{n-1}{3}},\enskip M^{(7)}(n)=\frac{169}{216}6^{\frac{n-1}{3}}.
\end{equation}
It is interesting that in this case $ps[\Lambda_9^3]$ is not used up to $M^{(7)},$ but the continuation of this list requires the knowing of $ps[\Lambda_{10}^3].$\newline
3) $j=2,\enskip n\geq32.$
$$M^{(1)}(n)=\frac{9}{4}6^{\frac{n-2}{3}},\enskip M^{(2)}(n)=\frac{13}{6}6^{\frac{n-2}{3}},$$

$$M^{(3)}(n)=2\cdot6^{\frac{n-2}{3}},\enskip M^{(4)}(n)=\frac{13}{9}6^{\frac{n-2}{3}},$$

$$M^{(5)}(n)=\frac{49}{36}6^{\frac{n-2}{3}},\enskip M^{(6)}(n)=\frac{4}{3}6^{\frac{n-2}{3}},$$

$$M^{(7)}(n)=\frac{31}{24}6^{\frac{n-2}{3}},\enskip M^{(8)}(n)=\frac{81}{64}6^{\frac{n-2}{3}},$$\newpage

\begin{equation}\label{8.9}
M^{(9)}(n)=\frac{5}{4}6^{\frac{n-2}{3}},\enskip M^{(10)}(n)=\frac{11}{9}6^{\frac{n-2}{3}},\enskip M^{(11)}(n)=\frac{39}{32}6^{\frac{n-1}{3}}.
\end{equation}
Note that the method not only gives a possibility to calculate the upper magnitudes of the permanent in $\Lambda_n^3,$ but also indicates those direct products on which they are attained. E.g., in (\ref{8.9}) $M_9$ is attained on direct products of some matrices $A_i\in \Lambda_i^3:$
 $$A_8\otimes \underbrace{A_3\otimes...\otimes A_3}_{\frac{n-8}{3}};\enskip A_4\otimes A_7\otimes \underbrace{A_3\otimes...\otimes A_3}_{\frac{n-11}{3}}; \enskip A_4\otimes A_4\otimes A_6\otimes \underbrace{A_3\otimes...\otimes A_3}_{\frac{n-14}{3}}.$$
Note also that the comparison of (\ref{8.7})-(\ref{8.9}) with (\ref{6.6})-(\ref{6.15}) shows that the following calculated $M^{(i)}$ are attained in $\widehat{\Lambda}_n^3,\enskip n\geq32: $\newline
in case $n\equiv0\mod{3},$ $$M^{(1)},M^{(2)},M^{(3)},M^{(4)};$$
in case $n\equiv1\mod{3},$
$$M^{(1)},M^{(3)},M^{(4)},M^{(5)},M^{(6)},M^{(7)}$$
(and is not attained $M^{(2)}$);\newline
in case $n\equiv2\mod{3},$
$$M^{(1)},M^{(2)},M^{(5)},M^{(7)},M^{(8)},M^{(9)},M^{(11)}$$
(and are not attained $M^{(3)},M^{(4)},M^{(6)},M^{(10)}$).
\section{Algorithm of a testing the parity of values of the permanent in $\Lambda_n^3$}

It seems that, among all known methods of calculation of the permanent, only Ryser's method (cf. \cite{11}, Ch.7) could be used for a creating an algorithm of a testing the parity of values of the permanent. Let $A$ be $n\times n$-matrix.  Let $A_r$ be a matrix which is obtained by changing some $r$ columns of $A$ by zero columns. Denote $S(A_r)$ the product of row sums of $A_r.$ Then, by Ryser's formula, we have
$$per A=\sum S(A_{0})-\sum S(A_{1})+$$
\begin{equation}\label{9.1}
\sum S(A_{2})-...+(-1)^{n-1}\sum S(A_{n-1}).
\end{equation}
Let now $A$ have integer elements. Introduce the following matrix function
\begin{equation}\label{9.2}
\Upsilon (A)=\begin{cases} 1,& if\enskip all \enskip row\enskip sums \enskip of \enskip A\enskip are\enskip odd,\\0, & otherwise.\end{cases}
\end{equation}
From (\ref{9.1}) we have
$$per A\equiv\sum \Upsilon(A_{0})-\sum \Upsilon(A_{1})+$$\newpage
\begin{equation}\label{9.3}
\sum \Upsilon(A_{2})+...+\sum \Upsilon(A_{n-1})\pmod {2}.
\end{equation}
Using (\ref{9.3}), let us create an algorithm of a search of the odd values of the permanent in $\Lambda_n^3.$ Since, evidently, $per A\equiv det A \pmod{2},$ then $A$ should have pairwise distinct columns. Note that cases $n\equiv j\pmod3, \enskip j=0,1,2,$ are considered by the same way. Suppose, say, $n=3t.$ According to (\ref{9.3}), we are interested in only cases when after removing $r\geq1$ columns of $A,$ all row sums will be odd. Suppose that after removing $r$ columns of $A,$ we have that $p$ sums remain to equal to 3 and $n-p$ sums equal to 1. This means that the total number of the removed 1's equals to $2(n-p)=6t-2p.$ Since, removing a column, we remove three 1's, then the number of the removed columns equals to $r=2t-\frac {2p}{3}.$ Thus $p=3m$ and $r=2(t-m), \enskip m=0,1,...,t-1.$
However, if $m=t-1,$ then $r=2.$ By the condition, these two columns are distinct, therefore, we conclude that at least one row sum equals to 2. The contradiction shows that the testing sequence is $r=4,6,...,2t.$ In cases $n\equiv1,2\pmod{3}$ we obtain the same testing sequence.
\begin{example}\label{e7}
Let us check the parities of values of the permanent of circulants in $\Delta_7^3\subset \Lambda_7^3.$
\end{example}
In this case $t=\lfloor\frac{7}{3}\rfloor=2$ and, therefore, the testing sequence contains only term $r=4.$
Note that matrix $A_r$ has all odd rows if and only if one row sum equals to 3 and each of 6 other row sums equals to 1. Indeed, let after the removing $4$ columns of $A,$ remain $p$ sums equal to 3 and $7-p$ sums equal to 1. This means that the total number of the removed 1's equals to $2(7-p)$ and the number of the removed columns equals to $r=4=\frac {14-2p} {3},$ i.e., $p=1.$ Moreover, since in a circulant all rows are congruent shifts of the first one, it is sufficient to consider the case when precisely the first row sum equals to 3 and others equal to 1 (the multiplication on 7 does not change the parity of the result). This opens a possibility of a momentary handy test on the parity every circulant of class $\Delta_7^3.$ This test consists of the removing all four columns beginning with 0. If now every rows $2,...,7$ has one 1, then the permanent is even; otherwise, it is odd. We check now directly that from $\binom{7}{3}=35$ circulants exactly 21 ones have odd permanent. \newline
\begin{remark}
\indent In $1967,$ Ryser \cite{14} did a conjecture that the number of the transversals of a latin square from elements $1,...,n$ ( i.e., the number of subsets of its $n$ pairwise distinct elements, none in the same row or column) has the same parity as $n.$ If $n$ is even, then the conjecture has been proved\newpage by Balasubramanian \cite{1}. Besides, in \cite{1} Balasubramanian did a conjecture for the parity of a sum of permanents, such that the truth of this conjecture yields Ryser's hypothesis for odd $n.$ In the same year $(1990),$ using the result of Example $7,$ the author disproved Balasubramanian's conjecture (private communication to Brualdi). It is interesting that soon Parker (see\cite{5}, p.$258$) indeed found several latin squares of order $7$ with even number of transversals. Add that later (\cite{18}) we found even an infinite set of counterexamples to the Balasubramanian conjecture.
\end{remark}
\section{Open problems}
1. To prove the MCI (Section 7).\newline
2. (Cf.\cite{17}, pp.171-172). Consider class $\Lambda_n(1,1+a,1+b),$ where $0\leq a\leq b<\frac {4e-9}{6}.$ Since
$\Lambda_n(1,1,1)=\Lambda_n^3,$ then  Voorhoeve's lower estimate for the permanent (\ref{2.16}) trivially holds for matrices in $\Lambda_n(1,1+a,1+b).$ It is clear that, for $a>0,\enskip b>0,$ it should exist an essentially stronger lower estimate. However, using Van der Waerden-Egorychev-Falikman theorem to class $\Lambda_n(\frac {1} {3+a+b},\frac{1+a}{3+a+b},\frac{1+b}{3+a+b})$ of doubly stochastic matrices, for the permanent of $\Lambda_n(1,1+a,1+b)$-matrices we obtain even weaker lower estimate of the order $C_1\sqrt{n}(\frac{3+a+b}{e})^n<<C(\frac{4}{3})^n.$ The problem is to find a stronger lower estimate for the permanent
in $\Lambda_n(1,1+a,1+b).$ \newline
3. (Cf.\cite{17}, pp.115-116). Let $M$ be a circulant of order n with integer elements. We conjecture that, for every integer $m,$ we have $per M\equiv(-1)^n per(mJ_n-M)\pmod{n},$ where $n\times n$-matrix $J_n$ consists of 1's only.
A special case of this conjecture, for $m=1,\enskip M=I_n+P+...+P^{k-1}$ in the equivalent terms was formulated by Yamamoto \cite{27} and proved for $k\leq3.$ The author \cite{15} proved the truth of the conjecture in case $m=1$ for arbitrary circulant $M$ ( including Yamamoto's conjecture for every $k).$ In \cite{17} the conjecture was proved for every $m$ and prime $n.$ The question is open in case of composite $n$ even in case $k=3.$\newline
4. Two Latin rectangles let us call equivalent, if the sets of their elements in the corresponding columns are the same. Note that numbers $|\Lambda_n^3|$ one can treat as the numbers of equivalence classes of Latin triangles. Let $A=I_n+P+P^2.$ In \cite{20} the author proved that the cardinality of the corresponding equivalent class is $2^n+6+2(-1)^n.$ To find the cardinality of the equivalent class which is defined by matrix $I_n+P+P^3.$

\end{document}